\documentclass[12pt]{amsart}
\usepackage{amssymb}
\usepackage{epsf}

\newcommand{\QED}{
\setlength{\unitlength}{1.0pt}%
\begin{picture}(7.5,7.5)(0,0)
\put(  0,  0){\line(1,0){7.5}}    \put(  0,2.5){\line(1,0){7.5}}
\put(  0,  5){\line(1,0){7.5}}    \put(2.5,7.5){\line(1,0){5}}
\put(  0,  0){\line(0,1){5}}    \put(2.5,  0){\line(0,1){7.5}}
\put(  5,  0){\line(0,1){7.5}}  \put(7.5,  0){\line(0,1){7.5}}
\put(  5,  0){\rule{2.5pt}{5pt}}
\end{picture}}

\newtheorem{thm}{Theorem}
\newtheorem{prop}[thm]{Proposition}
\newtheorem{cor}[thm]{Corollary}
\newtheorem{rem}[thm]{Remark}
\newtheorem{lemma}[thm]{Lemma}

\begin{document}

\title[A sagbi basis for the quantum Grassmannian]{A sagbi basis for the
quantum Grassmannian} 

\author{Frank Sottile}
\address{\hskip-\parindent
        Frank Sottile\\
	Department of Mathematics\\
        University of Wisconsin\\
        Madison, Wisconsin 53706\\
        USA}
\email{sottile@math.wisc.edu}

\author{Bernd Sturmfels}
\address{\hskip-\parindent
        Bernd Sturmfels\\
	Department of Mathematics\\
        University of California\\
        Berkeley, California 94720 \\
        USA}
\email{bernd@math.berkeley.edu}
\thanks{Second author supported in part by NSF grant DMS-9796181.
        Research at MSRI supported in part by NSF grant DMS-9701755}

\dedicatory{Dedicated to the memory of Gian-Carlo Rota}

\begin{abstract}
The maximal minors of a $p\times (m +  p)$-matrix of univariate
polynomials of degree $n$ with indeterminate coefficients
are themselves polynomials of degree $np$.
The subalgebra generated by their coefficients is
the coordinate ring of the quantum Grassmannian,
a singular compactification of the space of rational curves 
of degree $np$ in the Grassmannian of $p$-planes in ($m  +  p$)-space.
These subalgebra generators  are shown to form a sagbi basis.
The resulting flat deformation from the quantum Grassmannian  
to a toric variety gives a new ``Gr\"obner basis style'' proof of
the Ravi-Rosenthal-Wang formulas  in quantum Schubert calculus.
The coordinate ring of the quantum Grassmannian is an
algebra with straightening law, which is normal,
Cohen-Macaulay, Gorenstein and Koszul, 
and the ideal of quantum Pl\"ucker relations has a 
quadratic Gr\"obner basis.  This holds more
generally for skew quantum Schubert varieties. These results 
are well-known for the classical Schubert varieties $(n=0)$.
We also show that the row-consecutive $p\times p$-minors of
a generic matrix form a sagbi basis
and we give a quadratic Gr\"obner basis for their algebraic relations.
\end{abstract}
\maketitle

\section{Statement of the main result}
Let ${\mathcal M}(t)$ be the $p\times(m+p)$-matrix whose $i,j$th entry is
the degree $n$ polynomial in $t$,
$$
  x_{i,j}^{(n)}  \cdot t^n \ + \, 
  x_{i,j}^{(n-1)}  \cdot t^{n-1} +\,\,  \ldots \,\, + \,
  x_{i,j}^{(2)} \cdot t^2 \ + \,
  x_{i,j}^{(1)}  \cdot t\ + \,
  x_{i,j}^{(0)} . 
$$
The coefficients $x_{i,j}^{(l)}$ are indeterminates.
We write $k[X]$ for the polynomial ring over a field $k$
generated by these indeterminates, for 
$i=1,\ldots,p$, $j=1,\ldots,m+p$, and $l=0,\ldots,n$.
Lexicographic order on the triples $l,i,j$ 
gives a total order of these variables.
For example,
$$
  x_{1,2}^{(0)}\ <\ x_{1,2}^{(1)}\ <\ 
  x_{1,5}^{(1)}\ <\ x_{2,3}^{(1)}\ <\ 
  x_{2,4}^{(1)}\ <\ x_{1,3}^{(2)} \ .
$$
Let $\prec$ be the resulting degree reverse lexicographic term order on the 
polynomial ring $k[X]$.

\begin{figure}[htb]
 $$\epsfxsize=1.6in \epsfbox{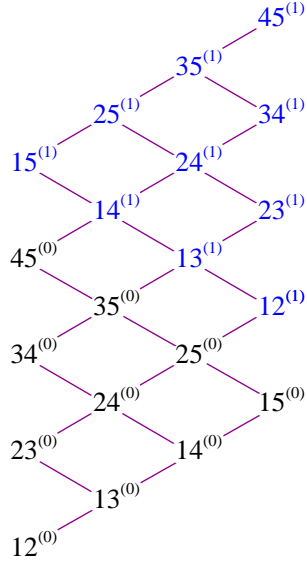}$$
 \caption{The distributive lattice ${\mathcal C}^1_{2,3}$\label{fig:qorder}.}
\end{figure}

For each $\alpha\in\binom{[m+p]}{p}$ and $a\geq 0$, let $\alpha^{(a)}$
be a variable which, when $a\leq np$,  formally represents the coefficient
of $t^a$ 
in the maximal minor of ${\mathcal M}(t)$ given by the columns indexed by
$\alpha_1,\alpha_2,\ldots,\alpha_p$. 
These variables $\alpha^{(a)}$ have a natural partial order,
denoted ${\mathcal C}_{p,m}$, which is defined as follows:
$$
  \alpha^{(a)}\ \leq\ \beta^{(b)} \quad\Longleftrightarrow\quad
   a\leq b \, \mbox{ and }  \,\,
   \alpha_i\leq \beta_{b-a+i} \,
  \mbox{ for } i = 1,2,\ldots,p-b+a. 
$$
Fix $0\leq q\leq np$ and let
${\mathcal C}^q_{p,m}$ denote the truncation of the infinite poset
${\mathcal C}_{p,m}$ to the finite subset $\,\bigl\{ \,
\alpha^{(a)} \,\, | \,\, \alpha\in\binom{[m+p]}{p} \mbox{ and } a\leq q 
\, \bigr\}$.
The posets ${\mathcal C}^q_{p,m}$ are graded distributive lattices.
Figure~\ref{fig:qorder} shows ${\mathcal C}^1_{2,3}$.

\noindent
Let $\varphi : k[{\mathcal C}^q_{p,m}] \rightarrow k[X]$
denote the $k$-algebra homomorphism
which sends the formal variable $\alpha^{(a)}$ to the coefficient of $t^a$
in the $\alpha$th maximal minor of the matrix ${\mathcal M}(t)$.
For example,
{\small 
$$ \begin{array}{l}
  \varphi \bigl( 456^{(2)} \bigr) \quad =
\quad \mbox{\rm coefficient of $t^2$ in}\  \
  \det\left[\begin{array}{ccc}
   x^{(0)}_{1,4}+x^{(1)}_{1,4}t&\ x^{(0)}_{1,5}+x^{(1)}_{1,5}t\ 
     & x^{(0)}_{1,6}+x^{(1)}_{1,6}t\\         \rule{0pt}{17pt}
   x^{(0)}_{2,4}+x^{(1)}_{2,4}t&\ x^{(0)}_{2,5}+x^{(1)}_{2,5}t\ 
     & x^{(0)}_{2,6}+x^{(1)}_{2,6}t\\         \rule{0pt}{17pt}
   x^{(0)}_{3,4}+x^{(1)}_{3,4}t&\ x^{(0)}_{3,5}+x^{(1)}_{3,5}t\ 
     & x^{(0)}_{3,6}+x^{(1)}_{3,6}t
  \end{array}\right]\,
\\ = \quad \rule{0pt}{20pt}%
  -\underline{x^{(0)}_{3,6}x^{(1)}_{1,5}x^{(1)}_{2,4}}
  +x^{(0)}_{3,5}x^{(1)}_{1,6}x^{(1)}_{2,4}
  +x^{(0)}_{3,6}x^{(1)}_{1,4}x^{(1)}_{2,5}
  -x^{(0)}_{3,4}x^{(1)}_{1,6}x^{(1)}_{2,5}
  -x^{(0)}_{3,5}x^{(1)}_{1,4}x^{(1)}_{2,6}
  +x^{(0)}_{3,4}x^{(1)}_{1,5}x^{(1)}_{2,6}\\\rule{0pt}{17pt}%
\phantom{=} \quad
  +x^{(0)}_{2,6}x^{(1)}_{1,5}x^{(1)}_{3,4}
  -x^{(0)}_{2,5}x^{(1)}_{1,6}x^{(1)}_{3,4}
  -x^{(0)}_{2,6}x^{(1)}_{1,4}x^{(1)}_{3,5}
  +x^{(0)}_{2,4}x^{(1)}_{1,6}x^{(1)}_{3,5}
  +x^{(0)}_{2,5}x^{(1)}_{1,4}x^{(1)}_{3,6}
  -x^{(0)}_{2,4}x^{(1)}_{1,5}x^{(1)}_{3,6}\\\rule{0pt}{17pt}%
\phantom{=} \quad
  -x^{(0)}_{1,6}x^{(1)}_{2,5}x^{(1)}_{3,4}
  +x^{(0)}_{1,5}x^{(1)}_{2,6}x^{(1)}_{3,4}
  +x^{(0)}_{1,6}x^{(1)}_{2,4}x^{(1)}_{3,5}
  -x^{(0)}_{1,4}x^{(1)}_{2,6}x^{(1)}_{3,5}
  -x^{(0)}_{1,5}x^{(1)}_{2,4}x^{(1)}_{3,6}
  +x^{(0)}_{1,4}x^{(1)}_{2,5}x^{(1)}_{3,6}\,.
 \end{array}
$$}

\begin{thm} \label{issagbi}
The set of polynomials $\, \varphi(\alpha^{(a)}) \,$ as 
$\alpha^{(a)} $ runs over the poset ${\mathcal C}^q_{p,m} \,$ forms a
sagbi basis with respect to the reverse lexicographic term order $\prec$
on $\,k[X]\,$ defined above.
\end{thm}

Our second theorem states that the subalgebra 
{\it image}$(\varphi)$ of $k[X]$ generated by
this sagbi basis is an {\it algebra with straightening law}
on the poset ${\mathcal C}^q_{p,m}$.
Let $\prec$ be the degree reverse lexicographic term order on 
$k[{\mathcal C}^q_{p,m}]$ induced by any linear extension of the poset 
${\mathcal C}^q_{p,m}$.
This term order on  $k[{\mathcal C}^q_{p,m}]$
and the previous term order on $k[X]$
are fixed throughout this paper.

\begin{thm}\label{thm:gbasis}
The reduced Gr\"obner basis of the kernel of $\varphi$
consists of quadratic polynomials in $\,k[{\mathcal C}^q_{p,m}]\,$ 
 which are indexed by 
pairs of incomparable variables $\gamma^{(c)},\delta^{(d)}$ in the poset 
$\,{\mathcal C}^{np}_{p,m}$,
$$
 S(\gamma^{(c)},\delta^{(d)}) \quad = \quad
  \gamma^{(c)}\cdot\delta^{(d)}\ -\ 
  (\gamma^{(c)}\vee\delta^{(d)})\cdot (\gamma^{(c)}\wedge\delta^{(d)})
\,\, + \,\,\hbox{lower terms in $\prec$},
$$
and all lower terms $\,\lambda\beta^{(b)}\alpha^{(a)}\,$ 
in $\,S(\gamma^{(c)},\delta^{(d)}) \,$
satisfy $\,\beta^{(b)}<\gamma^{(c)}\wedge\delta^{(d)}$ and 
$\gamma^{(c)}\vee\delta^{(d)}<\alpha^{(a)}$.
\end{thm}

The join $\vee$ and meet $\wedge$ appearing in the
above formula are the lattice operations 
in ${\mathcal C}^{np}_{p,m}$.
The combinatorial structure of this distributive lattice will become
clear in Section 2, when we introduce the toric variety associated
with  ${\mathcal C}^{np}_{p,m}$. In Section 3 we interpret the
subalgebra ${\rm image}(\varphi)$ of $k[X]$ as the
coordinate ring of the quantum Grassmannian. 
Section 4 contains the proofs of Theorems 1 and 2.
These results generalize the classical
sagbi basis property of maximal minors
\cite[Theorem 3.2.9]{Sturmfels_invariant}
and its geometric interpretation as a toric deformation
 \cite[Proposition 11.10]{Sturmfels_GBCP}
from the case of the Grassmannian to the quantum Grassmannian.
In Section 5 we discuss corollaries, applications
and some open problems.
One such application is that the row-consecutive $p\times p$-minors of any
matrix of indeterminates form a sagbi basis.\smallskip

We thank Aldo Conca, Ezra Miller, and Brian Taylor for their helpful
comments.

\section{The toric variety of the distributive lattice}

Theorem 1 asserts that the initial algebra 
of our subalgebra  $\,{\rm image}(\varphi)\,$
is generated by the initial monomials of its
generators $\varphi(\alpha^{(a)}) $.
Our first step is to identify the initial monomials.
Here are two examples. The first one is the underlined monomial
right before Theorem \ref{issagbi}:
$$
  \mbox{\rm in}_\prec \bigl( \varphi(456^{(2)}) \bigr) \ =\ 
  x_{3,6}^{(0)}x_{1,5}^{(1)}x_{2,4}^{(1)}\qquad\mbox{and}\qquad
  \mbox{\rm in}_\prec \bigl( \varphi(2457^{(5)}) \bigr) \ =\ 
  x_{2,7}^{(1)}x_{3,5}^{(1)}x_{4,4}^{(1)}x_{1,2}^{(2)}.
$$
In general, the initial monomial of $\varphi(\alpha^{(a)}) $
is given by the following lemma:

\begin{lemma}\label{lemma:leadmon}
Let $\alpha\in\binom{[m+p]}{p}$ and $\,a=pl+r \,$
with integers $p > r \geq 0$. Then
$$
  \mbox{\rm in}_\prec \bigl( \varphi(\alpha^{(a)}) \bigr) \quad =\quad 
  x_{r+1,\alpha_p}^{(l)}x_{r+2,\alpha_{p-1}}^{(l)}\cdots
  x_{p,\alpha_{r+1}}^{(l)}
x_{1,\alpha_r}^{(l+1)}
x_{2,\alpha_{r-1}}^{(l+1)}
\cdots
  x_{r,\alpha_1}^{(l+1)}.
$$
\end{lemma}

\noindent{\bf Proof. }
Let $x_{i_1,j_1}^{(l_1)}x_{i_2,j_2}^{(l_2)}\cdots x_{i_p,j_p}^{(l_p)}$
be a monomial which appears in $\varphi(\alpha^{(a)})$.
We claim that
\begin{equation}\label{eq:comparison}
  x_{i_1,j_1}^{(l_1)}x_{i_2,j_2}^{(l_2)}\cdots x_{i_p,j_p}^{(l_p)}
  \quad \preceq\quad
  x_{r+1,\alpha_p}^{(l)}x_{r+2,\alpha_{p-1}}^{(l)}\cdots
  x_{p,\alpha_{r+1}}^{(l)}x_{1,\alpha_r}^{(l+1)}\cdots
  x_{r,\alpha_1}^{(l+1)}.
\end{equation}
We may assume
$\,x_{i_1,j_1}^{(l_1)} \prec x_{i_2,j_2}^{(l_2)} \prec \cdots \prec
x_{i_p,j_p}^{(l_p)}\,$ and hence $l_1 \leq\cdots\leq l_p$.
Since $l_1+\cdots+l_p=a$, either 
$ l_1 < q $,
from which~(\ref{eq:comparison}) follows, or else 
$l_1=\cdots=l_{p-r}=l$ and $l_{p+1-r}=\cdots=l_p=l+1$.
In the second case, as $\{i_1,\ldots,i_p\}=\{1,\ldots,p\}$ and the monomial
is in order, we must have $i_1<\cdots<i_{p-r}$ and $i_{p+1-r}<\cdots<i_p$.
If  $i_1 \leq r$,
then~(\ref{eq:comparison}) follows, and 
if $ i_1 = r+1$,
then the ordered sequence
$\,i_1,i_2,\ldots,i_p\,$
 equals $\, r \! + \!  1,r \! + \! 2,\ldots,p,1,\ldots,r$.
Among all monomials satisfying this new second case, the largest in the
degree reverse lexicographic order $\prec $ has 
the second lower index appearing in reverse order.
This completes the proof.
\QED\medskip

We next introduce some combinatorics to help understand the poset
${\mathcal C}_{p,m}$.
A {\it row} with shift $a$ consists of $p$ consecutive empty unit boxes
shifted $a$ units to the left of a given vertical line.
A {\it skew shape} is an array of such rows whose shifts are weakly
increasing read from top to bottom.
For example, the unshaded boxes in the figure on the left form a skew
shape (with shifts 1,2,2, and 5) while those in the other figure do not:
$$
    \epsfxsize=2.8in \epsfbox{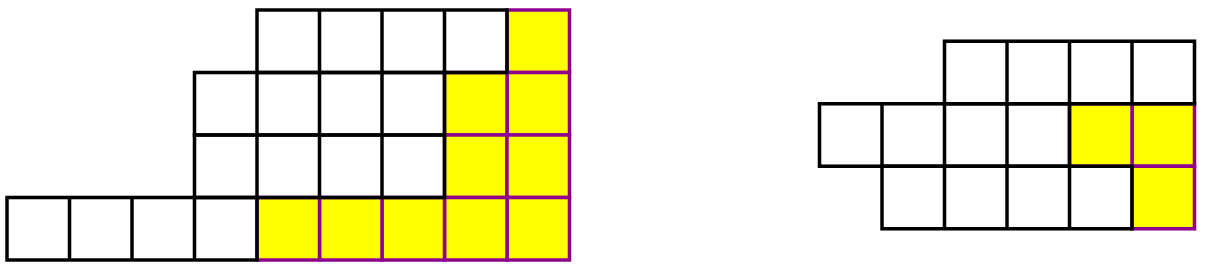}
$$

A {\it (skew) tableau} $T$ is a filling of a skew shape with integers that
increase across each row.
When the entries lie in $[m+p]$, the $i$th row of a tableau 
is a sequence $\alpha_{(i)}\in\binom{[m+p]}{p}$.
If $a_i$ is the shift of the  $i$th row and $T$ has $j$ rows, then 
$T$ corresponds to a monomial 
$\alpha_{(1)}^{(a_1)}\alpha_{(2)}^{(a_2)}\cdots\alpha_{(j)}^{(a_j)}$.
Conversely, any monomial in the variables $\alpha^{(a)}$ corresponds to a
tableau.
A tableau $T$ is {\it standard} if the entries are weakly increasing in each
column, read top to bottom.
Equivalently, $T$ is standard if we have
$\alpha_{(1)}^{(a_1)}\leq\alpha_{(2)}^{(a_2)}\leq\cdots
        \leq\alpha_{(j)}^{(a_j)}$
in ${\mathcal C}_{p,m}$.
For example, the following two tableaux correspond to the monomials
$345^{(0)}123^{(1)}245^{(3)}$ and $135^{(0)}123^{(1)}257^{(3)}$.
The first tableau is not standard and the second tableau is standard:
$$
    \epsfxsize=2.8in \epsfbox{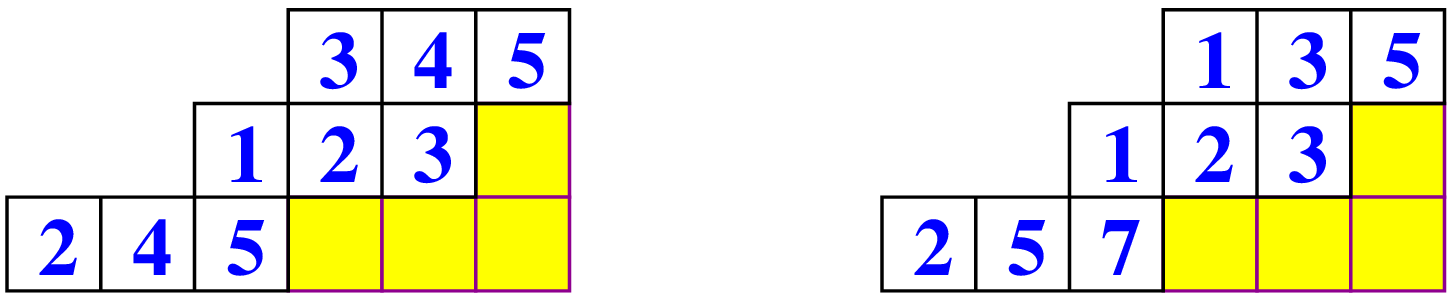}
$$

The elements of the poset ${\mathcal C}^q_{p,m}$
are represented by one-row tableaux with entries in 
$[m+p]$ and shift at most $q$.
Two elements satisfy $\,\alpha^{(a)} \leq \beta^{(b)} \,$
if and only if the two-rowed tableau
$T=\alpha^{(a)}\beta^{(b)}$ is standard.
This representation implies that ${\mathcal C}^q_{p,m}$
is a distributive lattice. 
Indeed, the two lattice operations
$\wedge$ and $\vee$ are described as follows.
If a two-rowed tableau $\,T=\alpha^{(a)}\beta^{(b)}\,$ 
is non-standard then
interchanging the entries in every column in which a violation 
($\alpha_{a-b+i}<\beta_i$) occurs yields a standard tableau.
The first row of this new tableau is the {\it meet} 
$\alpha^{(a)}\wedge\beta^{(b)}$ of  
$\alpha^{(a)}$ and $\beta^{(b)}$ in ${\mathcal C}^q_{p,m}$ and the 
second row is their {\it join} $\alpha^{(a)}\vee\beta^{(b)}$.

Let $\psi : k[{\mathcal C}^{q}_{p,m}] \rightarrow k[X]$
denote the $k$-algebra homomorphism
which sends the variable $\alpha^{(a)}$ to the
monomial $\,\mbox{\rm in}_\prec \bigl(\varphi(\alpha^{(a)}) \bigr)$.
Its kernel is a {\it toric ideal} \
(i.e.~binomial prime) \ in $\, k[{\mathcal C}^{q}_{p,m}]$.

\begin{prop}\label{lem:toric}
The reduced Gr\"obner basis for the kernel of 
$\psi$ consists of the binomials
$$
\underline{\alpha^{(a)}\cdot\beta^{(b)}} \ -\ 
  (\alpha^{(a)}\vee\beta^{(b)})\cdot(\alpha^{(a)}\wedge\beta^{(b)}),
$$
where $\alpha^{(a)}, \beta^{(b)}$ runs over all incomparable pairs of \/
${\mathcal C}^q_{p,m} $. The initial monomial is underlined.
\end{prop}

\noindent{\bf Proof.}
This follows from Hibi's Theorem \cite{Hibi} since 
${\mathcal C}^q_{p,m}$ is a distributive lattice.
\QED

\begin{cor}\label{cor:standard}
The set of standard tableaux is a $k$-basis for 
$\, k[{\mathcal C}^q_{p,m}] / {\rm kernel}(\psi) = 
{\rm image}(\psi) $.
\end{cor}

Here is a typical element in the reduced Gr\"obner basis of
${\rm kernel}(\psi) $  for $p=5,m=4,q=9$:
$$
\underline{45789^{(1)} \cdot 12356^{(3)}} \, - \,
35689^{(1)} \cdot 12457^{(3)}.$$
Note that the second monomial corresponds
to a standard tableau while the first does not.

We write $T^q_{p,m}$ for the  projective toric
variety cut out by the binomials in Lemma \ref{lem:toric}.
Its coordinate ring is the subalgebra $\,{\rm image}(\psi)$ of $k[X]$.
The geometry of toric varieties associated
with distributive lattices is discussed in
\cite{Wagner}. The analogue to Corollary \ref{cor:standard} always holds,
i.e., multichains in the poset correspond to basis 
monomials in the coordinate ring. 

\begin{cor} \label{Degree} 
The degree of the toric variety $T^q_{p,m}$ is
the number of maximal chains in ${\mathcal C}^q_{p,m}$.
\end{cor}

The closed intervals of the poset ${\mathcal C}^q_{p,m}$
are also distributive lattices. They are denoted
$$ [\beta^{(b)}, \alpha^{(a)}] \quad := \quad
\bigl\{ \gamma^{(c)} \in {\mathcal C}^q_{p,m} \,: \,
\beta^{(b)} \leq \gamma^{(c)} \leq  \alpha^{(a)} \bigr\}. $$
Proposition \ref{lem:toric} and Corollaries 
\ref{cor:standard} and  \ref{Degree} hold 
essentially verbatim for the distributive sublattice
$\,[\beta^{(b)}, \alpha^{(a)}]\,$ as well.
The projective toric variety associated with $\,[\beta^{(b)}, \alpha^{(a)}]\,$
is gotten from the toric variety of ${\mathcal C}^q_{p,m}$
by setting $\gamma^{(c)} = 0$ for all
$ \gamma^{(c)} \not\in [\beta^{(b)}, \alpha^{(a)}]$.
The degree of that variety is the number of 
saturated chains in ${\mathcal C}^q_{p,m}$ which 
start at $\beta^{(b)}$ and end at $ \alpha^{(a)}$.

We close this section with an alternative proof, to be used
in Section 4, for the
fact that ${\mathcal C}_{p,m}$ is a distributive lattice.
We claim that ${\mathcal C}_{p,m}$ is a sublattice
of {\it Young's lattice}. 
Given  $\alpha^{(a)} \in {\mathcal C}_{p,m}$, write
$a=pl+r$ with integers $p > r \geq 0$, and define
a sequence $J(\alpha^{(a)})$ by
\begin{equation}\label{eq:seq_def}
  J(\alpha^{(a)})_i\quad :=\quad \left\{\begin{array}{ll}
   l(m+p)+\alpha_{r+i}&\quad\mbox{if }1 \leq i\leq p-r\\
   (l+1)(m+p)+\alpha_{i-p+r}&\quad\mbox{if }p-r<i\leq p
  \end{array}\right.\,.
\end{equation}
This gives an order-preserving bijection between 
the poset ${\mathcal C}_{p,m}$ and
the poset of sequences $J:=j_1<j_2<\cdots<j_p$ of positive integers with 
$j_p-(m+p)<j_1$, and it preserves meet and join.
This bijection preserves the rank function in the two 
distributive lattices:
\begin{equation}\label{Miracle}
| \alpha^{(a)} | \quad := \quad
a(m+p) + \sum_{j=1}^p (\alpha_j - j) 
\quad = \quad
\sum_{i=1}^p \bigl(  J(\alpha^{(a)})_i - i \,\bigr)
\quad =: \quad
| J(\alpha^{(a)})|.
\end{equation}

\section{The quantum Grassmannian}

Let ${\it Grass}_pk^{m+p}$ denote the Grassmannian of
$p$-planes in the vector space $k^{m+p}$. This is a smooth
projective variety of dimension $mp$. 
Consider the space $S^q_{p,m}$ of maps 
${\mathbb P}^1\rightarrow \mbox{\it Grass}_pk^{m+p}$ of degree $q$.
Such a map may be (non-uniquely) represented as the row space of a
$p\times(m+p)$-matrix of polynomials in $t$ whose maximal minors have
degree $q$. Results in \cite{Clark} imply that it suffices to
consider the matrices ${\mathcal M}(t)$ in the introduction.
The coefficients of these maximal minors define the {\it Pl\"ucker
embedding} of $S^q_{p,m}$ into 
${\mathbb P}(\wedge^pk^{m+p}\otimes k^{q+1})$; see \cite{Stromme,Rosen94}.
The {\it quantum Grassmannian} $K^q_{p,m}$ is the Zariski closure of
$S^q_{p,m}$ in this Pl\"ucker embedding. It is an irreducible
projective variety of dimension $mp+q(m+p)$.  Its prime ideal
is $\,{\rm kernel}(\varphi) \subset k[ {\mathcal C}^q_{p,m}]\,$ and 
its coordinate ring is our subalgebra
$\,{\rm image}(\varphi)\subset k[X]$.

The quantum Grassmannian $K^q_{p,m}$
is singular and it differs from other spaces used to
study rational curves in Grassmann varieties (the quot
scheme~\cite{Stromme}, the Kontsevich space of stable
maps~\cite{Kontsevich_Manin}, or the set of autoregressive
systems~\cite{RR94}). 
Nevertheless, $K^q_{p,m}$ has been crucial in two 
important advances: in computing the intersection number 
$\,{\rm degree}(K^q_{p,m})\,$ in
quantum cohomology~\cite{RRW98}, and in showing that this
intersection problem can be fully solved over the real
numbers~\cite{Sottile_quantum}. 
Our result will give a new derivation of this intersection number.

\begin{cor}\label{cor:counting}
{\rm \cite{RRW98}}
The degree of $\,K^q_{p,m} \,$ is the number of maximal chains in
${\mathcal C}^q_{p,m}$. 
\end{cor}

\noindent{\bf Proof. }
This follows immediately from Theorem 2 and 
Corollary \ref{Degree}.
\QED\medskip

This degree can also be computed in the small quantum cohomology ring of the
Grassmannian~\cite{Bertram}.  Note that
$\, {\rm degree} (K^1_{2,3}) = 55 $, by
counting maximal chains in Figure~\ref{fig:qorder}.
Ravi, Rosenthal, and Wang \cite{RRW98}
were motivated by a problem in applied mathematics.
The degree of $\,K^q_{p,m} \,$ is the number of dynamic 
feedback compensators that stabilize a certain linear system,
in the sense of systems theory.
This number can be described in classical projective geometry
as follows.
The Schubert subvariety of ${\it Grass}_pk^{m+p}$ consisting of 
$p$-planes meeting a fixed $m$-plane $L$ is a hyperplane section in the
Pl\"ucker embedding of ${\it Grass}_pk^{m+p}$.
Thus the set of maps $M\in S^q_{p,m}$ such that 
$M(t)$ meets $L$ non-trivially is
a hyperplane section of $S^q_{p,m}$ in its Pl\"ucker embedding.
Since $GL_{m+p}(k)$ acts transitively on 
${\it Grass}_pk^{m+p}$, Kleiman's Theorem on generic
transversality  implies the following
statement when $k$ is algebraically closed of characteristic zero.
Set $N:=mp+q(m+p)$ and suppose  
$t_1,\ldots,t_N\in{\mathbb P}^1$ are general points and 
$L_1,\ldots,L_N$ are general $m$-planes in $k^{m+p}$, 
then the degree of $K^q_{p,m}$ counts those maps 
$M$ for which $M(t_i)$ meets $L_i$ non-trivially for each $i=1,\ldots,N$.
As to computing the  desired maps $M$ numerically, we
note that the sagbi basis in Theorem 1 and
the Gr\"obner basis in Theorem 2 each lead to an
{\it optimal homotopy algorithm} for finding these 
$\,{\rm degree}(K^q_{p,m})\,$ 
maps. These algorithms generalize the ones in~\cite{HSS}.

\begin{rem}\label{rem:flat}
The sagbi basis in Theorem 1 defines a flat
deformation from the quantum Grassmannian
$K^q_{p,m}$ to the projective toric variety
$T^q_{p,m}$ associated with the poset
${\mathcal C}^q_{p,m}$.
\end{rem}

See \cite{CHV} for a precise algebraic discussion of such deformations,
and see \cite[Equation (11.9)]{Sturmfels_GBCP} for the simplest
example relevant to us, namely, $K^0_{2,3} = {\it Grass}_2 k^{5}$.
The flat deformation is given algebraically
by deleting all but the first two terms in the
Gr\"obner basis elements  $\, S(\gamma^{(c)},\delta^{(d)}) \,$ 
given in Theorem 2.  Consider 
the deformation from $K^3_{3,3}$ to $T^3_{3,3}$.
The incomparable pair $156^{(1)}$ and $234^{(2)}$ in ${\mathcal C}^3_{3,3}\,$
indexes the quadratic
polynomial $\,\, S(156^{(1)},234^{(2)})\, = $
{\small 
\begin{equation}
\label{BigSyzygy}
\begin{array}{c}
\underline{ 156^{(1)}234^{(2)}
-146^{(1)}235^{(2)} }\,
+145^{(1)}236^{(2)}
+136^{(1)}245^{(2)}
-135^{(1)}246^{(2)}\\ 
+134^{(1)}256^{(2)}
-126^{(1)}345^{(2)} 
+125^{(1)}346^{(2)}
-124^{(1)}356^{(2)} 
+123^{(1)}456^{(2)}\\
-456^{(0)}123^{(3)} 
+356^{(0)}124^{(3)}
-346^{(0)}125^{(3)} 
+345^{(0)}126^{(3)}
-256^{(0)}134^{(3)} \\ 
+246^{(0)}135^{(3)}
-245^{(0)}136^{(3)} 
-236^{(0)}145^{(3)}
+235^{(0)}146^{(3)}
-234^{(0)}156^{(3)} \\
+2{\cdot}156^{(0)}234^{(3)}
-2{\cdot}146^{(0)}235^{(3)} 
+2{\cdot}145^{(0)}236^{(3)}
+2{\cdot}136^{(0)}245^{(3)}
-2{\cdot}135^{(0)}246^{(3)}\,\\
+2{\cdot}134^{(0)}256^{(3)}
-2{\cdot}126^{(0)}345^{(3)} 
+2{\cdot}125^{(0)}346^{(3)}
-2{\cdot}124^{(0)}356^{(3)} 
+2{\cdot}123^{(0)}456^{(3)},
\end{array}
\end{equation}
}\vspace{-6pt}

\noindent
which vanishes on the quantum Grassmannian $K^3_{3,3}$.
The underlined leading binomial vanishes on the toric variety $T^3_{3,3}$,
by Proposition \ref{lem:toric}.
Our main technical problem, to be solved in the next section,
is the reconstruction of quadrics such as (\ref{BigSyzygy})
from their leading binomial.

A key tool in proving Theorems 1 and 2
is the  Schubert decomposition 
of the quantum Grassmannian $K^q_{p,m}$
indexed by ${\mathcal C}^q_{p,m}$.
For $\alpha^{(a)}\in{\mathcal C}^q_{p,m}$, the {\it quantum Schubert
variety} is
$$
  Z_{\alpha^{(a)}}\quad :=\quad 
\bigl\{ \,({\gamma^{(c)}})\in K^q_{p,m} \,\mid\,
{\gamma^{(c)}}=0 \, \mbox{ if } \, \gamma^{(c)}\not\leq \alpha^{(a)} \bigr\}\ .
$$
More generally, 
for $\beta^{(b)}\leq \alpha^{(a)}$ in 
${\mathcal C}^q_{p,m}$, we define the {\it skew quantum Schubert variety}
$$
  Z_{\alpha^{(a)}/\beta^{(b)}}\quad :=\quad 
    \bigl\{\, ({\gamma^{(c)}})\in K^q_{p,m} \,\mid \,
{\gamma^{(c)}}  = 0 \, \mbox{ if } \,
      \gamma^{(c)}\not\in [\beta^{(b)},\alpha^{(a)}] \bigr\}\ .
$$
Among the quantum Schubert varieties of $K^q_{p,m}$ are the 
$K^d_{p,m}$ for $d<q$; namely, if
 $\delta^{(d)}$ is the supremum of ${\mathcal C}^d_{p,m}$, then 
$K^d_{p,m}=Z_{\delta^{(d)}}$.
This allows us to deduce assertions about the general
quantum Grassmannian $K^q_{p,m}$ from
results about quantum Schubert varieties of $K^{pn}_{p,m}$.

The quantum Schubert varieties and skew quantum Schubert varieties
have rational parameterizations which are constructed as follows.
Let $\alpha^{(a)}\in{\mathcal C}^{pn}_{p,m}$ and write $a=ps+r$ with 
integers $p>r\geq 0$.   We define the matrix
 ${\mathcal M}_{\alpha^{(a)}}(t)$ to be the specialization of 
${\mathcal M}(t)$ where
$$
  x_{i,j}^{(l)}\ =\ 0\quad\mbox{if }\quad\left\{
  \begin{array}{ll}
   ( l>s+1 \  \mbox{ and } \  i\leq r )
&\mbox{ or } \ \ (l=s+1\mbox{ and }j>\alpha_{r+1-i}) \ \ \mbox{or} \\
( l>s  \ \mbox{ and } \ i>r )
&\mbox{ or } \ \ (l=s\mbox{ and }j>\alpha_{p+r+1-i})\,.
  \end{array}\right.
$$
Here we use the conventions
$\alpha_\nu = 0$ if $\nu \leq 0\, $ and
$\,\alpha_\nu = +\infty$ if $\nu > p$.
For example,
$$
{\mathcal M}_{235^{(2)}}(t) \quad = \quad \,
 \left[\begin{array}{cccccc}
  x_{1,1}^{(0)}+x_{1,1}^{(1)}\cdot t&x_{1,2}^{(0)}+x_{1,2}^{(1)}\cdot t
  &x_{1,3}^{(0)}+x_{1,3}^{(1)}\cdot t&x_{1,4}^{(0)}
  &\ x_{1,5}^{(0)}&\ x_{1,6}^{(0)}\\ \rule{0pt}{17pt}
  x_{2,1}^{(0)}+x_{2,1}^{(1)}\cdot t&x_{2,2}^{(0)}+x_{2,2}^{(1)}\cdot t 
     &x_{2,3}^{(0)}\ \ \ &\ x_{2,4}^{(0)}\ 
     &\ x_{2,5}^{(0)} &\ x_{2,6}^{(0)} \\
  \rule{0pt}{17pt}
  x_{3,1}^{(0)}\ \ \ &x_{3,2}^{(0)}\ \ \ 
     & x_{3,3}^{(0)}\ \ \ &\ x_{3,4}^{(0)}\ &\ x_{3,5}^{(0)} &\ \ 0
  \rule{0pt}{16pt}
 \end{array}\right].
$$
If we specialize the variables $x^{(l)}_{i,j}$ in 
${\mathcal M}_{\alpha^{(a)}}(t)$ to field elements in $k$ 
in such a way that the  resulting matrix over $k(t)$
has maximal row rank, then that matrix defines a map
from $k$ to ${\it Grass}_pk^{m+p}$.
If we extend this to ${\mathbb P}^1$, we obtain a map in 
$Z_{\alpha^{(a)}}$. 
Proposition~\ref{ItIsDominant} below implies that 
such maps constitute a dense subset of $Z_{\alpha^{(a)}}$.
This means that the coefficients with respect to $t$ of the 
maximal minors of ${\mathcal M}_{\alpha^{(a)}}(t)$ give a rational
parameterization of $Z_{\alpha^{(a)}}$.  

This construction extends to skew quantum Schubert varieties as follows.
Given $\beta^{(b)}\leq\alpha^{(a)}$, write $b=ps+r$ with integers
$p>r\geq 0$  and define the  matrix 
${\mathcal M}_{\alpha^{(a)}/\beta^{(b)}}(t)$ to be the specialization of 
${\mathcal M}_{\alpha^{(a)}}(t)$ where
$$
  x_{i,j}^{(l)}\ =\ 0\quad\mbox{if }\quad\left\{
  \begin{array}{ll}
( l<s+1   \  \mbox{ and }\  i\leq r )
    &\mbox{ or } \ \ 
   ( l=s+1\mbox{ and }j< \beta_{r+1-i})
\ \  \mbox{ or }  \\
( l<s   \ \mbox{ and }\  i>r ) 
    &   \mbox{ or } \ \ 
   ( l=s\mbox{ and }j<\beta_{p+r+1-i})
  \end{array}\right.\,.
$$
The matrix ${\mathcal M}_{\alpha^{(a)}/\beta^{(b)}}(t)$ gives a rational
map into $Z_{\alpha^{(a)}/\beta^{(b)}}$, which is
described algebraically as follows.
 We define $\varphi_{\alpha^{(a)}}$ and $\varphi_{\alpha^{(a)}/\beta^{(b)}}$ 
to be the composition of the map
$\,\varphi : k[{\mathcal C}^{np}_{p,m}] \rightarrow k[X] \,$
with the specializations to
${\mathcal M}_{\alpha^{(a)}}(t)$ and 
${\mathcal M}_{\alpha^{(a)}/\beta^{(b)}}(t)$ respectively.
We claim that these matrices parameterize dense subsets of the 
(skew) quantum Schubert varieties.

\begin{prop} \label{ItIsDominant}
The kernel of $\varphi_{\alpha^{(a)}}$ is the homogeneous ideal of the
quantum Schubert variety $Z_{\alpha^{(a)}}$.
Likewise, the  kernel of $\varphi_{\alpha^{(a)}/\beta^{(b)}}$ is the
homogeneous ideal of the skew quantum Schubert variety
$Z_{\alpha^{(a)}/\beta^{(b)}}$. In particular, the varieties 
$Z_{\alpha^{(a)}}$ and $Z_{\alpha^{(a)}/\beta^{(b)}}$ are irreducible.
\end{prop}

We postpone the proof of this proposition until the next section.
Here is an example which illustrates the parameterization of skew 
quantum Schubert varieties for $p = m = 3$:
$$
{\mathcal M}_{235^{(2)}/146^{(1)}}(t) \quad = \quad
 \left[\begin{array}{cccccc}
  x_{1,1}^{(1)}\cdot t&x_{1,2}^{(1)}\cdot t
  &x_{1,3}^{(1)}\cdot t&0&0&0\\ \rule{0pt}{17pt}
  x_{2,1}^{(1)}\cdot t&x_{2,2}^{(1)}\cdot t 
     &0&0&0 &x_{2,6}^{(0)} \\
  \rule{0pt}{17pt}
  0&0&0&x_{3,4}^{(0)}&x_{3,5}^{(0)}& 0
  \rule{0pt}{16pt}
 \end{array}\right]\,.
$$
We evaluate the $3 \times 3$-minors of this matrix to find the
$k$-algebra homomorphism $\,\varphi_{{235^{(2)}/146^{(1)}}} $. It 
takes polynomials in $12$ variables
$\,\gamma^{(c)}\,$ to polynomials
in $8$ variables $x^{(l)}_{i,j}$ as follows:
$$ 
 \begin{array}{c}
  146^{(1)} \mapsto  -x^{(1)}_{1,1} x^{(0)}_{2,6} x^{(0)}_{3,4} \,\,, \quad
  156^{(1)} \mapsto  -x^{(1)}_{1,1} x^{(0)}_{2,6} x^{(0)}_{3,5} \,\, ,\quad 
  246^{(1)} \mapsto  -x^{(1)}_{1,2} x^{(0)}_{2,6} x^{(0)}_{3,4} \,,
\\\rule{0pt}{15pt}
  256^{(1)} \mapsto  -x^{(1)}_{1,2} x^{(0)}_{2,6} x^{(0)}_{3,5} \,\, ,\quad 
  346^{(1)} \mapsto  -x^{(1)}_{1,3} x^{(0)}_{2,6} x^{(0)}_{3,4} \,\, ,\quad 
  356^{(1)} \mapsto  -x^{(1)}_{1,3} x^{(0)}_{2,6} x^{(0)}_{3,5} \, , 
\\\rule{0pt}{15pt}
  124^{(2)} \mapsto  x^{(1)}_{1,1} x^{(1)}_{2,2} x^{(0)}_{3,4}
                   - x^{(1)}_{2,1} x^{(1)}_{1,2} x^{(0)}_{3,4} \, , \quad
  125^{(2)} \mapsto  x^{(1)}_{1,1} x^{(1)}_{2,2} x^{(0)}_{3,5}
                   - x^{(1)}_{2,1} x^{(1)}_{1,2} x^{(0)}_{3,5}\, , 
\\\rule{0pt}{15pt}
  134^{(2)} \mapsto  -x^{(1)}_{2,1} x^{(1)}_{1,3} x^{(0)}_{3,4} \,\, ,\qquad 
  135^{(2)} \mapsto  -x^{(1)}_{2,1} x^{(1)}_{1,3} x^{(0)}_{3,5} \,, 
\\\rule{0pt}{15pt}
  234^{(2)} \mapsto  -x^{(1)}_{2,2} x^{(1)}_{1,3} x^{(0)}_{3,4} \,\, ,\qquad 
  235^{(2)} \mapsto  -x^{(1)}_{2,2} x^{(1)}_{1,3} x^{(0)}_{3,5} .
\end{array}\vspace{-4pt}
$$
The $12$ variables $\gamma^{(c)}$ appearing on the
left sides above are precisely the elements in the interval
$\,\bigl[146^{(1)},  235^{(2)} \bigr]\,$
of the distributive lattice  ${\mathcal C}_{3,3}$. 
There are $18$ incomparable pairs in this interval,
each giving a quadratic generator for the
kernel of $\,\varphi_{{235^{(2)}/146^{(1)}}} $. 
This set of $18$ quadrics consists of $14$ binomials
and four trinomials, and it equals the reduced
Gr\"obner basis with respect to $\prec$.
For example, one of the 14 binomials in this Gr\"obner basis 
is the underlined leading binomial of 
$\,\, S(156^{(1)},234^{(2)})\, $ in~(\ref{BigSyzygy}), and 
one of the four trinomials is
$$ \underline{ 346^{(1)} \cdot  125^{(2)} } \,\, - \,\,
 246^{(1)} \cdot  135^{(2)} \,\, + \,\,
 146^{(1)} \cdot  235^{(2)}.
$$
The underlined term is an incomparable pair in
$\bigl[146^{(1)},  235^{(2)} \bigr]$, while the other two
monomials are comparable pairs. Erasing the third term
gives a  binomial as  in Proposition \ref{lem:toric}.

\section{Construction of Straightening Syzygies}\label{sec:syzygy}

The following theorem is the technical heart of this paper.
All three of Theorem \ref{issagbi}, Theorem \ref{thm:gbasis}, and
Proposition \ref{ItIsDominant} will be derived from
Theorem  \ref{thm:syzygy} in the end of this section.

\begin{thm}\label{thm:syzygy}
Let $\gamma^{(c)},\delta^{(d)}$ be a pair of incomparable variables 
in the poset ${\mathcal C}^{np}_{p,m}$.
There is a  quadric $S(\gamma^{(c)},\delta^{(d)})$ in the kernel
of $\varphi : k[{\mathcal C}^{np}_{p,m}] \rightarrow k[X]$
 whose first two monomials are
$$
  \gamma^{(c)}\cdot\delta^{(d)}\ -\ 
  (\gamma^{(c)}\vee\delta^{(d)}) \cdot (\gamma^{(c)}\wedge\delta^{(d)}).
$$
Moreover, if $\lambda\beta^{(b)}\alpha^{(a)}$ 
is any non-initial monomial in $S(\gamma^{(c)},\delta^{(d)})$,
then
$\gamma^{(c)},\delta^{(d)}\in[\beta^{(b)},\alpha^{(a)}]$. 
\end{thm}

The pair $ \beta^{(b)}\alpha^{(a)}$  
in the second assertion is necessarily standard, 
i.e.~$\beta^{(b)}<\alpha^{(a)}$.
The quadrics $S(\gamma^{(c)},\delta^{(d)})$ 
are not constructed explicitly, but rather through an iterative
procedure modeled on the {\it subduction algorithm} in {\rm image}$(\varphi)$.
A main idea is to utilize the well-known subduction process
\cite[Algorithm 3.2.6]{Sturmfels_invariant}
modulo the $p \times p$-minors of a generic $p \times N$-matrix.

Set $N:=(n+1)(m+p)$.
Let ${\mathcal N}$ be the $p\times N$-matrix whose $i,j$th
entry is $x^{(l)}_{i,r}$, where 
$j= (m+p)l+r$ with $1\leq r\leq m+p$.
If ${\mathcal N}_l$ is the submatrix of ${\mathcal N}$ consisting of the
entries $x^{(l)}_{i,j}$, then ${\mathcal N}$ is the concatenation of
${\mathcal N}_0,{\mathcal N}_1,\ldots,{\mathcal N}_n$ and 
${\mathcal M}(t)=
{\mathcal N}_0+t{\mathcal N}_1+\cdots+t^n{\mathcal N}_n$.
Sequences $\,J:j_1<\cdots<j_p\in\binom{[ N]}{p}\,$
are regarded as variables. We write
$\phi(J)$ for the $J$th maximal minor of ${\mathcal N}$.
{\it Young's poset} on sequences $J$ is given by
componentwise comparison and is graded via
$\,|J| := \sum_i (j_i-i)$. 
The coefficient $\varphi(\alpha^{(a)})$ of $t^a$ in the $\alpha$th maximal
minor of ${\mathcal M}(t)$ is an
alternating sum of maximal minors of ${\mathcal N}$.
The exact formula is
\begin{equation}\label{eq:phi_expand}
  \varphi(\alpha^{(a)})\quad \,\,\, =\quad 
  \sum_{\stackrel{\mbox{\scriptsize $|J|=|\alpha^{(a)}|$}}%
      {J\equiv\alpha\,\bmod{(m+p)}}}
\!\!  \epsilon_J \cdot \phi(J),
\end{equation}
where $\epsilon_J$ is the sign of the permutation
that orders the following sequence:
$$
  j_1\bmod(m+p),\ j_2\bmod(m+p),\ \ldots,\ j_p\bmod(m+p).
$$

The polynomial 
rings $k[{\mathcal C}^{np}_{p,m}]$ and $k[\binom{[N]}{p}]$ are graded
with $\deg\alpha^{(a)}=|\alpha^{(a)}|$ and $\deg J=|J|$.
Consider the degree-preserving $k$-algebra homomorphism
$\pi:k[{\mathcal C}^{np}_{p,m}]\rightarrow k[\binom{[N]}{p}]$ defined by
\begin{equation}\label{eq:pidef}
  \pi(\alpha^{(a)})\quad \,\, \, =\quad 
  \sum_{\stackrel{\mbox{\scriptsize $|J|=|\alpha^{(a)}|$}}%
      {J\equiv\alpha\,\bmod{(m+p)}}}
\!\!  \epsilon_J \cdot J\,.
\end{equation}
Lexicographic order on the sequences $J\in\binom{[ N]}{p}$
gives a linear extension of Young's poset.
In this ordering, the initial term of~(\ref{eq:pidef}) is
the sequence $J(\alpha^{(a)})$ defined in (\ref{eq:seq_def}).
This sequence is characterized by 
$\, {\rm in}_\prec \,\varphi(\alpha^{(a)}) \, = \,{\rm in}_\prec
\, \phi(J(\alpha^{(a)}))$.
It can be checked that all other terms 
$ \epsilon_J \cdot J \,$ appearing in~(\ref{eq:pidef})
satisfy $\, J_1 < J(\alpha^{(a)})_1 \,$
and $\,  J_p-J_1  \ >\ m+p $.
For example, for
$m = 4$,
$$
\pi(235^{(2)}) \quad = \quad
\underline{   (5, 9, 10)} 
-   (3, 9, 12)
+   (3, 5, 16)
+   (2, 10, 12)
- (2, 5, 17)
+ (2,3,19),
$$
$$  \mbox{and} \qquad
{\rm in}_\prec \bigl( \varphi( \,235^{(2)} \,) \bigr) \quad = \quad
{\rm in}_\prec\bigl(\phi( 5,9,10) \bigr) \quad = \quad
x_{3,5}^{(0)}
x_{1,3}^{(1)}
x_{2,2}^{(1)} . \qquad
 $$

For $J\in\binom{[ N]}{p}$, let ${\mathcal N}_J$ be the
specialization of ${\mathcal N}$ where in each row $i$, all entries in
columns greater than $j_i$ are set to zero.
Under the identification of ${\mathcal N}$ with ${\mathcal M}(t)$, we have
${\mathcal N}_{J(\alpha^{(a)})}={\mathcal M}_{\alpha^{(a)}}$.
Let $\phi_J : k [\binom{[N]}{p}] \rightarrow k[X] $ denote the
$k$-algebra homomorphism which maps the formal variable $I$ to the
$I$th maximal minor of ${\mathcal N}_J$.
Then $\phi_J(I)$ vanishes unless $I\leq J$.
In particular, if $|I|=|J|$, then $\phi_J(I)$ vanishes unless $I=J$, and in
that case, it is just the product of the last non-zero variables in each row
of ${\mathcal N}_J$. 
{}From this it follows that
\begin{equation}
\label{SomeIdentities}
  \begin{array}{rcl}
& \varphi_{\alpha^{(a)}} \quad = \quad 
\phi_{J(\alpha^{(a)})} \circ \pi  \\\rule{0pt}{15pt}
& \varphi_{\alpha^{(a)}}(\alpha^{(a)}) \,\, =\,\,
   \phi_{J(\alpha^{(a)})}(J(\alpha^{(a)})) \,\, = \,\,
   {\rm in}_\prec \, \varphi(\alpha^{(a)}) \,\, =
\,\, \psi(\alpha^{(a)}). 
\end{array}
\end{equation}
In the Pl\"ucker embedding of {\it Grass}$_pk^{ N}$ into 
${\mathbb P}(\wedge^p k^{N})$, the Schubert variety indexed by $J$ 
is 
$$
  \Omega_J\quad :=\quad \{y=(y_I)\in {\it Grass}_pk^{ N}\mid
  y_I=0\mbox{ if }I\not\leq J \}\ .
$$
The homogeneous ideal ${\mathcal I}(\Omega_J)$  
which defines this Schubert variety
is precisely the kernel of $\phi_J$.
The following identity of ideals in $k[\binom{[N]}{p}]$
follows from the classical Pl\"ucker relations:

\begin{prop}\label{prop:Schubert-ideal}
For any $J\in\binom{[ N]}{p}$ we have
$$
  \bigcap_{I<J}{\mathcal I}(\Omega_I)\quad =\quad 
  {\mathcal I}(\Omega_J) \, + \, \left\langle\, J \,\right\rangle\,.
$$
\end{prop}

The map $\pi:k[{\mathcal C}^{np}_{p,m}]\rightarrow k[\binom{[N]}{p}]$
induces a birational isomorphism 
$\pi^*:{\it Grass}_pk^{ N}\dashrightarrow K^{np}_{p,m}$.
{}From the identification of ${\mathcal M}_{\alpha^{(a)}}(t)$ with 
${\mathcal N}_{J(\alpha^{(a)})}$ and Proposition~\ref{ItIsDominant},
we will see that 
$\pi^*(\Omega_{J(\alpha^{(a)})})$ is a dense subset of $\, Z_{\alpha^{(a)}}$.
We also consider the image under $\pi^*$ 
of the Schubert varieties $\Omega_J$ for $\, J<J(\alpha^{(a)})$.

\begin{prop}[Ravi-Rosenthal-Wang~\cite{RRW98}]\label{prop:O_J-image}
If $J<J(\alpha^{(a)})$, then 
$$   \pi^*(\Omega_J)\quad \subset\quad 
  \bigcup_{\beta^{(b)}<\alpha^{(a)}} Z_{\beta^{(b)}}\,.
$$
\end{prop}

\noindent{\bf Proof. }
The inclusion $\Omega_J\subset\Omega_{J(\alpha^{(a)})}$ implies
$\pi^*(\Omega_J)\subset Z_{\alpha^{(a)}}$.
Since $\varphi_{\alpha^{(a)}}(\alpha^{(a)})$ is the product of leading
entries in the rows of ${\mathcal N}_{J(\alpha^{(a)})}$, it follows that 
$\varphi_{\alpha^{(a)}}(\alpha^{(a)})$ vanishes under the specialization 
to ${\mathcal N}_J$, and hence $\pi(\alpha^{(a)})$ vanishes on
$\Omega_J$. 
This implies our claim because
$\bigcup_{\beta^{(b)}<\alpha^{(a)}} Z_{\beta^{(b)}}$
is defined as a subvariety of $ Z_{\beta^{(b)}}$ by the vanishing of
$\alpha^{(a)}$.
\QED\medskip

For $L<J$ in Young's poset, 
define ${\mathcal N}_{J/L}$ to be the specialization of
${\mathcal N}$ where in the $i$th row, only the entries in columns
$l_i,l_i+1,\ldots,j_i$ are non-zero.
Then ${\mathcal M}_{\alpha^{(a)}/\beta^{(b)}}(t)$ is the
specialization of ${\mathcal M}(t)$ corresponding to 
${\mathcal N}_{J(\alpha^{(a)})/J(\beta^{(b)})}(t)$.
Define the $k$-algebra homomorphism
$\,\phi_{J/I} : k [\binom{[N]}{p}] \rightarrow k[X] \,$ 
by evaluating the appropriate minors on
${\mathcal N}_{J/L}$.
We observe that
\begin{equation}\label{eq:initial}
  \begin{array}{rcl}
    \varphi_{\alpha^{(a)}/\beta^{(b)}}(\alpha^{(a)}) &=&
     {\rm in}_\prec \, \varphi(\alpha^{(a)}) \quad =
 \quad \phi_{J(\alpha^{(a)})/J(\beta^{(b)})}\bigl( J(\alpha^{(a)})\bigr),
\\\rule{0pt}{16pt}
    \varphi_{\alpha^{(a)}/\beta^{(b)}}(\beta^{(b)}) &=&
     {\rm in}_\prec  \, \varphi(\beta^{(b)})
\quad = \quad \phi_{J(\alpha^{(a)})/J(\beta^{(b)})}\bigl( J(\beta^{(b)})\bigr).
  \end{array}
\end{equation}

The following lemma is very useful in our proof of Theorem~\ref{thm:syzygy}.  

\begin{lemma}\label{lem:factor}
Fix $\alpha^{(a)} \in {\mathcal C}^{np}_{p,m}$ and
let $f\in k[{\mathcal C}^{np}_{p,m}]$ be a quadratic form of degree $d$.
\begin{enumerate}
\item
     Suppose that $\varphi_{\beta^{(b)}}(f)=0$ for all
     $\beta^{(b)}<\alpha^{(a)}$.
     Then there exist constants $\lambda_J\in k$ with 
   $$
     \varphi_{\alpha^{(a)}}(f)\quad =\quad 
     \varphi_{\alpha^{(a)}}(\alpha^{(a)}) \,\, \cdot  \!\!\!\!
     \sum_{\stackrel{\mbox{\scriptsize $J \in \binom{[ N]}{p}$}}%
      {|J|+|\alpha^{(a)}|=d}}
       \lambda_J \cdot \phi_{J(\alpha^{(a)})}(J)\,.
   $$
\item
     Suppose $\beta^{(b)}<\alpha^{(a)}$ and 
     $\varphi_{\alpha^{(a)}/\gamma^{(c)}}(f)=0$ for all
     $\beta^{(b)}<\gamma^{(c)}\leq\alpha^{(a)}$.
    For some $\lambda_J\in k$,
   $$
     \varphi_{\alpha^{(a)}/\beta^{(b)}}(f)\quad =\quad 
     \varphi_{\alpha^{(a)}/\beta^{(b)}}(\beta^{(b)}) \, \cdot 
\!\!\!\!      \sum_{\stackrel{\mbox{\scriptsize $J \in \binom{[ N]}{p}$}}%
      {|J|+|\beta^{(b)}|=d}} \lambda_J\cdot 
     \phi_{J(\alpha^{(a)})/J(\beta^{(b)})}(J)\,.
   $$
\end{enumerate}
\end{lemma}

\noindent{\bf Proof. }
We only prove part 1. The hypothesis states  that
$\, \phi_{J(\alpha^{(a)})}(\pi(f)) \, = \, 
\varphi_{\alpha^{(a)}}(f) \, $
vanishes on all matrices 
${\mathcal N}_{J(\beta^{(b)})}$ for
     $\beta^{(b)}<\alpha^{(a)}$.
Proposition \ref{prop:O_J-image} implies that
$\pi(f)$ vanishes on all 
Schubert varieties $\Omega_J$ with $J<J(\alpha^{(a)})$.
But then, using Proposition
\ref{prop:Schubert-ideal},
$$
  \pi(f)\quad \in\quad 
  \bigcap_{J<J(\alpha^{(a)})}{\mathcal I}(\Omega_J)\quad =\quad 
  {\mathcal I}(\Omega_{J(\alpha^{(a)})}) \, + \,
  \left\langle \,     J(\alpha^{(a)}) \, \right\rangle.
$$
This means $\,\pi(f)=g + J(\alpha^{(a)})\cdot h$, 
where $g\in {\mathcal I}(\Omega_{J(\alpha^{(a)})})=
{\rm ker} \bigl(\phi_{J(\alpha^{(a)})} \bigr) $
and $h\in k[\binom{[N]}{p}]$ is a linear form of degree $d-|\alpha^{(a)}|$.
Such a linear form can be written as follows
$$
  h \, \quad =\quad \sum_{|J|+|\alpha^{(a)}|=d} \! \! \lambda_J\,J\,.
$$
By applying the map $\, \phi_{J(\alpha^{(a)})} \,$ to both
sides of the equation $\,\pi(f)=g + J(\alpha^{(a)})\cdot h$, 
we obtain the first assertion of Lemma \ref{lem:factor}.
Part 2 is proved by similar arguments.
\QED\medskip

Our proof of Theorem~\ref{thm:syzygy} will show that the sums in
Lemma \ref{lem:factor} are actually sums of terms of the form
$\,\lambda_{J(\delta^{(d)})}\cdot\varphi_{\alpha^{(a)}}(\delta^{(d)})\,$ 
and  $\,\lambda_{J(\delta^{(d)})}\cdot
        \varphi_{\alpha^{(a)}/\beta^{(b)}}(\delta^{(d)}) \,$
respectively. The next lemma provides the initial step in our inductive 
proof of Theorem~\ref{thm:syzygy}.

\begin{lemma}\label{lem:initial}
Let $\gamma^{(c)}$ and $\delta^{(d)}$ be incomparable variables in the poset 
${\mathcal C}^{np}_{p,m}$ and set
$\alpha^{(a)}:=\gamma^{(c)}\vee\delta^{(d)}$.
Then $
\varphi_{\alpha^{(a)}}(\gamma^{(c)}\cdot\delta^{(d)}\ -\  
  \gamma^{(c)}\vee\delta^{(d)}\cdot\gamma^{(c)}\wedge\delta^{(d)})
   =0$.
\end{lemma}

\noindent{\bf Proof. }
We prove the lemma by inductively showing that,
for each  $\,\beta^{(b)}\leq\alpha^{(a)}$,
\begin{equation}\label{eq:lt}
  \varphi_{\alpha^{(a)}/\beta^{(b)}}
\bigl( \, \gamma^{(c)}\cdot\delta^{(d)}\ -\  
  \gamma^{(c)}\vee\delta^{(d)}\cdot\gamma^{(c)}\wedge\delta^{(d)} \,\bigr)
   \quad  = \quad 0.
\end{equation}
If $\beta^{(b)}\not\leq\gamma^{(c)}\wedge\delta^{(d)}$, then 
$\varphi_{\alpha^{(a)}/\beta^{(b)}}(\gamma^{(c)}\wedge\delta^{(d)})$
vanishes, and either
$\varphi_{\alpha^{(a)}/\beta^{(b)}}(\gamma^{(c)})$ vanishes or
$\varphi_{\alpha^{(a)}/\beta^{(b)}}(\delta^{(d)})$ vanishes.
This implies that~(\ref{eq:lt}) holds. 

Next suppose $\beta^{(b)}=\gamma^{(c)}\wedge\delta^{(d)}$. We claim that  
$\varphi_{\alpha^{(a)}/\beta^{(b)}}$ maps each variable appearing
in~(\ref{eq:lt}) to its initial term in $k[X]$. In view of
Proposition~\ref{lem:toric}, this claim implies~(\ref{eq:lt}).
To establish this claim, we need only show that
$\varphi_{\alpha^{(a)}/\beta^{(b)}}(\gamma^{(c)})=
  {\rm in}_\prec\varphi(\gamma^{(c)})$, as the case for $\delta^{(d)}$ is
similar and that of the other terms follow from~(\ref{eq:initial}).
Consider the expansion of 
$\varphi_{\alpha^{(a)}/\beta^{(b)}}(\gamma^{(c)})$ in terms of the minors
$\phi(J)$ of 
${\mathcal N}_{J(\alpha^{(a)})/J(\beta^{(b)})}$.
First observe that the submatrix given by the columns from $J(\gamma^{(c)})$
is block anti-diagonal, with each block either upper or lower
triangular along its anti-diagonal.
This is because for each $i$, $J(\gamma^{(c)})_i$ is either 
$J(\beta^{(b)})_i$ or $J(\alpha^{(a)})_i$, and the non-zero entries in the
$i$th row of  
${\mathcal N}_{J(\alpha^{(a)})/J(\beta^{(b)})}$ lie between these two
numbers.
Thus the contribution of term $J(\gamma^{(c)})$ to 
$\varphi_{\alpha^{(a)}/\beta^{(b)}}(\gamma^{(c)})$ is simply 
${\rm in}_\prec \, \varphi(\gamma^{(c)})$.

We claim there are no other terms.
If there were another term indexed by $L$, then the $L$th maximal minor of 
${\mathcal N}_{J(\alpha^{(a)})/J(\beta^{(b)})}$ would be non-zero, and so
$J(\beta^{(b)})\leq L\leq J(\alpha^{(a)})$.
Thus
$$
  J(\beta^{(b)})\quad =\quad J(\gamma^{(c)})\wedge J(\delta^{(d)})
  \quad \leq\quad L\wedge J(\delta^{(d)}).
$$
Comparing the first components of these sequences gives 
$\min\{J(\gamma^{(c)})_1,J(\delta^{(d)})_1\}\leq L_1$.
Since $L_1<J(\gamma^{(c)})_1$, this implies 
$J(\delta^{(d)})_1\leq L_1$.
Similarly, using $J(\alpha^{(a)})\geq L\vee J(\delta^{(d)})$, we see that
$J(\delta^{(d)})_p\geq L_p$.
Lastly, as $L$ is a summand in $\pi(\gamma^{(c)})$ and 
$L\neq J(\gamma^{(c)})$, we have $L_p-L_1> m+p$ and thus 
$$
  m+p\ \geq\ J(\delta^{(d)})_p-J(\delta^{(d)})_1\ \geq\ 
  L_p-L_1\ >\ m+p,
$$
a contradiction, which proves the claim.
Thus~(\ref{eq:lt}) holds for $\beta^{(b)}=\gamma^{(c)}\wedge\delta^{(d)}$.

Finally, let $\zeta^{(z)}<\gamma^{(c)}\wedge\delta^{(d)}$ and
suppose that~(\ref{eq:lt}) holds for all $\beta^{(b)}$ with 
$\zeta^{(z)}<\beta^{(b)}\leq \alpha^{(a)}$.
Then by Lemma~\ref{lem:factor},
$$
  \varphi_{\alpha^{(a)}/\zeta^{(z)}}\bigl(\gamma^{(c)}\cdot\delta^{(d)}\ -\  
    \gamma^{(c)}\vee\delta^{(d)}\cdot\gamma^{(c)}\wedge\delta^{(d)} \bigr)
   \quad =\quad 
  \varphi_{\alpha^{(a)}/\zeta^{(z)}}(\zeta^{(z)})\cdot
  \sum_{J}\lambda_J \cdot
     \phi_{J(\alpha^{(a)})/J(\zeta^{(z)})}(J)\,,
$$
the sum over sequences $J$ of rank 
$|J|=|\gamma^{(c)}|+|\delta^{(d)}|-|\zeta^{(z)}|$.
But this exceeds the rank of $\alpha^{(a)}$, since
$\zeta^{(z)}<\gamma^{(c)}\wedge\delta^{(d)}$ and  
$|\alpha^{(a)}|+|\gamma^{(c)}\wedge\delta^{(d)}|
=|\gamma^{(c)}|+|\delta^{(d)}|$. 
Thus the sum vanishes and so~(\ref{eq:lt}) holds for all
$\beta^{(b)}\leq\alpha^{(a)}$, which proves the lemma.
\QED\medskip

\noindent{\bf Proof of Theorem~\ref{thm:syzygy}. }
Let $\gamma^{(c)}$ and $\delta^{(d)}$ be incomparable variables in the poset 
${\mathcal C}^{np}_{p,m}$.
For each $\alpha^{(a)}\in{\mathcal C}^{np}_{p,m}$ we inductively construct
quadratic polynomials 
$S_{\alpha^{(a)}}(\gamma^{(c)},\delta^{(d)})\in
k[\beta^{(s)}\mid\beta^{(s)}\leq \alpha^{(a)}]$, and then show 
$S_{\alpha^{(a)}}(\gamma^{(c)},\delta^{(d)})$ is in the kernel of the map
$\varphi_{\alpha^{(a)}}$. 
The case when $\alpha^{(a)}$ is the top element in the poset 
${\mathcal C}^{np}_{p,m}$ proves the theorem.
These polynomials have the following restriction property:
If $\beta^{(b)}<\alpha^{(a)}$, then 
$S_{\beta^{(b)}}(\gamma^{(c)},\delta^{(d)})$ is the image of
$S_{\alpha^{(a)}}(\gamma^{(c)},\delta^{(d)})$ under the map which sets
variables $\zeta^{(z)}\not\leq \beta^{(b)}$ to zero.
They also have the further homogeneity properties that each non-zero term 
$\lambda\zeta^{(z)}\beta^{(b)}$ must have
$z+b=c+d$ and satisfy the multiset equality
$\beta\cup\zeta=\gamma\cup\delta$, 
and if it is not the initial term, then 
$\gamma^{(c)},\delta^{(d)}\in[\zeta^{(z)},\beta^{(b)}]$.

For $\alpha^{(a)}\not\geq \gamma^{(c)}\vee\delta^{(d)}$, set
$S_{\alpha^{(a)}}(\gamma^{(c)},\delta^{(d)}):=0$
and if $\alpha^{(a)}=\gamma^{(c)}\vee\delta^{(d)}$, then set
$$
  S_{\alpha^{(a)}}(\gamma^{(c)},\delta^{(d)})\quad:=\quad
  \gamma^{(c)}\cdot\delta^{(d)}\ -\  
    \gamma^{(c)}\vee\delta^{(d)}\cdot\gamma^{(c)}\wedge\delta^{(d)}.
$$
These polynomials have the restriction and homogeneity properties, and, for
$\alpha^{(a)}\not> \gamma^{(c)}\vee\delta^{(d)}$, we have
$\varphi_{\alpha^{(a)}}(S_{\alpha^{(a)}}(\gamma^{(c)},\delta^{(d)}))=0$, by
Lemma~\ref{lem:initial}. 

Let $\alpha^{(a)}>\gamma^{(c)}\vee\delta^{(d)}$ and suppose we have
constructed $S_{\beta^{(b)}}(\gamma^{(c)},\delta^{(d)})$ for each 
$\beta^{(b)}<\alpha^{(a)}$.
By the restriction property, there is a polynomial 
$S'\in k[\beta^{(b)}\mid \beta^{(b)}<\alpha^{(a)}]$ which restricts to 
$S_{\beta^{(b)}}(\gamma^{(c)},\delta^{(d)})$ for each 
$\beta^{(b)}<\alpha^{(a)}$.
Thus $\varphi_{\beta^{(b)}}(S')=0$ for all $\beta^{(b)}<\alpha^{(a)}$.

Set $e:=|\gamma^{(c)}|+|\delta^{(d)}|$, the degree of 
$S'$.
By Lemma~\ref{lem:factor}, 
\begin{equation}\label{eq:sprime}
  \varphi_{\alpha^{(a)}}(S')\quad =\quad 
  \varphi_{\alpha^{(a)}}(\alpha^{(a)})\cdot 
  \sum_{|J|+|\alpha^{(a)}|=e}\lambda_J \cdot\phi_{J(\alpha^{(a)})}(J)\,.
\end{equation}
If we consider the columns of ${\mathcal M}_{\alpha^{(a)}}(t)$ involved in
$\varphi_{\alpha^{(a)}}(S')$, we see that this 
sum is further restricted to those $J$ which satisfy the multiset equality
$(\gamma \cup\delta)\setminus \alpha\equiv J\mod(m+p)$, 
with $J\bmod(m+p)$ consisting of distinct integers, and with
$J\leq J(\alpha^{(a)})$.
If there are no such $J$, then 
$\varphi_{\alpha^{(a)}}(S')=0$ and we set 
$S_{\alpha^{(a)}}(\gamma^{(c)},\delta^{(d)})=S'$.

Otherwise, let $z:=c+d-a$ and $\zeta:=(\gamma\cup\delta)\setminus\alpha$.
Then the summands in~(\ref{eq:sprime}) are among those $J$ which
appear in $\pi(\zeta^{(z)})$ so we have $J(\zeta^{(z)})<J(\alpha^{(a)})$
and hence $\zeta^{(z)}<\alpha^{(a)}$.
Observe that $\varphi_{\alpha^{(a)}/\zeta^{(z)}}(S')=
  \lambda_{J(\zeta^{(z)})} \cdot
  \varphi_{\alpha^{(a)}/\zeta^{(z)}}(\alpha^{(a)}\zeta^{(z)})$.
Define
$$
  S_{\alpha^{(a)}}(\gamma^{(c)},\delta^{(d)})\quad :=\quad 
  S' - \lambda_{J(\zeta^{(z)})} \alpha^{(a)}\zeta^{(z)}.
$$
We claim that if $\lambda_{J(\zeta^{(z)})}\neq 0$, then 
$\zeta^{(z)}\leq\gamma^{(c)},\delta^{(d)}$.
If not, then every term of $S'$ contains a
variable $\xi^{(x)}$ with $\zeta^{(z)}\not\leq\xi^{(x)}$, and so we must
have $\varphi_{\alpha^{(a)}/\zeta^{(z)}}(S')=0$,
a contradiction.

We complete the proof of Theorem~\ref{thm:syzygy} by showing that for
$\beta^{(b)}\leq \alpha^{(a)}$, 
\begin{equation}\label{eq:indstep}
  \varphi_{\alpha^{(a)}/\beta^{(b)}}(S_{\alpha^{(a)}}
  (\gamma^{(c)},\delta^{(d)}))\quad =\quad 0.
\end{equation}
If $\beta^{(b)}\not\leq \zeta^{(z)}$, then
$\varphi_{\alpha^{(a)}/\beta^{(b)}}(\zeta^{(z)})=0$ and so
$\varphi_{\alpha^{(a)}/\beta^{(b)}}(S_{\alpha^{(a)}}
    (\gamma^{(c)},\delta^{(d)}))=\varphi_{\alpha^{(a)}/\beta^{(b)}}(S')$,
which is zero as 
$\phi_{J(\alpha^{(a)})/J(\beta^{(b)})}(J)=0$ for all $J$ which appear in 
$\pi(\zeta^{(z)})$.
By the construction of $S_{\alpha^{(a)}}(\gamma^{(c)},\delta^{(d)}))$, we
also have $\varphi_{\alpha^{(a)}/\zeta^{(z)}}(S_{\alpha^{(a)}}
    (\gamma^{(c)},\delta^{(d)}))=0$.

Let $\xi^{(x)}<\zeta^{(z)}$ and suppose~(\ref{eq:indstep}) holds for all 
$\beta^{(b)}$ with $\xi^{(x)}<\beta^{(b)}$.
Then by Lemma~\ref{lem:factor},
$$
  \varphi_{\alpha^{(a)}/\xi^{(x)}}(S_{\alpha^{(a)}}
    (\gamma^{(c)},\delta^{(d)}))\quad =\quad 
  \varphi_{\alpha^{(a)}/\xi^{(x)}}(\xi^{(x)})\cdot 
  \sum_{|J|+|\xi^{(x)}|=e} \lambda_J \cdot 
\phi_{J(\alpha^{(a)})/J(\xi^{(x)})}(J).
$$
Since $|J|=e-|\xi^{(x)}|>e-|\zeta^{(z)}|=|\alpha^{(a)}|$, 
each term in the right hand sum is zero. 
\QED\medskip\smallskip

It is now straightforward
to derive all our assertions that were left unproven so far.

\medskip

\noindent{\bf Proof of Theorem~\ref{issagbi}.}
Theorem  \ref{thm:syzygy} together with
Proposition \ref{lem:toric} shows that the 
subduction criterion for sagbi bases
(see e.g.~\cite[Proposition 1.1]{CHV} or
\cite[Theorem 11.4]{Sturmfels_GBCP}) is satisfied.
\QED  \medskip

\noindent{\bf Proof of Theorem~\ref{thm:gbasis}. }
A standard fact on sagbi bases, proved in
\cite[Corollary 2.2]{CHV} or in
\cite[Corollary 11.6 (1)]{Sturmfels_GBCP}, states that
the reduced Gr\"obner basis for the binomial ideal ${\rm kernel}(\psi)$
lifts to a reduced Gr\"obner basis for 
the non-binomial ideal ${\rm kernel}(\varphi)$.
\QED \medskip

\noindent{\bf Proof of Proposition \ref{ItIsDominant}. }
If $\beta^{(b)}$ is the minimal element in the poset
${\mathcal C}_{p,m}$, then $\varphi_{\alpha^{(a)}}$ and 
$\varphi_{\alpha^{(a)}/\beta^{(b)}}$ have the same kernel,
as the varieties $Z_{\alpha^{(a)}/\beta^{(b)}}$ and $Z_{\alpha^{(a)}}$ are
equal.  Hence  it suffices to prove the second statement
about $\varphi_{\alpha^{(a)}/\beta^{(b)}}$.
Clearly, the kernel of $\varphi_{\alpha^{(a)}/\beta^{(b)}}$ 
contains the homogeneous ideal of the skew quantum Schubert variety
$Z_{\alpha^{(a)}/\beta^{(b)}}$. 
If this containment were proper,
then we would also get proper containment at the level of 
initial ideals with respect to the induced partial term order,
which was denoted by ${\mathcal A}^T \omega$ in
\cite[Chapter 11]{Sturmfels_GBCP}. But that is impossible
since every binomial relation on the monomials
$\,{\rm in}_\prec \, \varphi_{\alpha^{(a)}/\beta^{(b)}}( \gamma^{(c)}) =
{\rm in}_\prec \, \varphi ( \gamma^{(c)}) \,$
lifts to a polynomial which vanishes on
$Z_{\alpha^{(a)}/\beta^{(b)}}$, as shown in the proof
Theorem  \ref{thm:syzygy}.
\QED

\section{Applications and future directions}

We first summarize some algebraic consequences of our main results.

\begin{cor}\label{ASL}
The coordinate ring of the quantum Grassmannian 
$K^q_{p,m}$ is an algebra with straightening law
on the distributive lattice $C^q_{p,m}$. 
It has a presentation
by a non-commutative Gr\"obner basis consisting 
of quadratic elements, and, in particular, it is a Koszul algebra.
\end{cor}

\noindent{\bf Proof.}
The first statement follows from Theorem~\ref{issagbi} and the form of the
syzygies $S(\gamma^{(c)},\delta^{(d)})$ of Theorem~\ref{thm:gbasis}.
For the second statement, 
consider the coordinate ring of $K^q_{p,m}$
as the quotient of the free associative algebra
on $C^q_{p,m}$ modulo a two-sided ideal. 
By \cite[Proposition 3.2]{EPS} that two-sided ideal has a quadratic 
Gr\"obner basis, obtained from lifting the Gr\"obner basis in 
Theorem~\ref{thm:gbasis}.
For the classical Grassmannian $(n=0)$ this result appeared in
\cite{Gr_Hu}. The Koszul property is a well-known consequence of the
existence of a quadratic Gr\"obner basis;
see e.g.~\cite[Theorem 3]{Gr_Hu} for the non-commutative 
version which is relevant here.
\QED

\begin{cor}\label{C-MG}
The coordinate ring of $K^q_{p,m}$ is
a normal Cohen-Macaulay and Gorenstein domain.
It has rational singularities if ${\rm char}(k) = 0$
and it is $F$-rational if ${\rm char}(k) > 0$.
\end{cor}

\noindent{\bf Proof. }
By Corollary~\ref{ASL} and the results of~\cite{CHV}, these properties of
$K^q_{p,m}$ follow from the corresponding properties of the toric variety
$T^q_{p,m}$.
But these were established in~\cite{Wagner}, as $T^q_{p,m}$ is the toric
variety associated to the distributive lattice ${\mathcal C}^q_{p,m}$.
\QED\medskip

We remark that both Corollary~\ref{C-MG} and the analog of
Corollary~\ref{ASL} (with the poset ${\mathcal C}^q_{p,m}$
replaced by the appropriate interval) hold for the skew quantum Schubert
varieties.

Our next application is the sagbi property of the
row-consecutive $p \times p$-minors of a matrix of indeterminates.
This result is non-trivial since the set of {\it all} $p\times p$-minors
is not a sagbi basis in general~\cite[Example 11.3]{Sturmfels_GBCP}.
A finite sagbi basis for the algebra
of all $p \times p$-minors was found
by Bruns and Conca \cite{Bruns_Conca}.
Let ${\mathcal L}$ be the $p(n+1)\times(m+p)$-matrix whose $i,j$th entry is
$x^{(l)}_{r,j}$, where $i=pl+r$.
This matrix is obtained from ${\mathcal N}$ by stacking the matrices
${\mathcal N}_0, \ldots, {\mathcal N}_n$.
Let $\chi:k[{\mathcal C}^{np}_{p,m}]\rightarrow k[X]$ denote the $k$-algebra
homomorphism which sends the variable $\alpha^{(a)}$ to the $\alpha$th
maximal minor of the submatrix of ${\mathcal L}$ consisting of rows 
$a+1,a+2,\ldots,a+p$.
Thus the collection of polynomials $\chi(\alpha^{(a)})$ are the
row-consecutive $p\times p$-minors of ${\mathcal L}$.

\begin{thm}\label{row-c-sagbi}
The set 
$\, \bigl\{ \chi(\alpha^{(a)}) \, : \, \alpha^{(a)} \in
{\mathcal C}^{np}_{p,m} \bigr\} \,$
of row-consecutive $p \times p$-minors of a generic matrix
is a sagbi basis with respect 
to the degree reverse lexicographic term order $\prec$ on $k[X]$.
\end{thm}

Brain Taylor has pointed out this may also be deduced from
Proposition~2.7.3 of his Ph.D.~Thesis~\cite{BDT_thesis}.\medskip

\noindent{\bf Proof. }
Let $\omega$ be the weight on the variables in $k[X]$  defined by 
$\omega(x^{(l)}_{i,j}):= -(pl+i)^2$.
Then $\chi(\alpha^{(a)}) = {\rm in}_\omega
\bigl( \varphi(\alpha^{(a)})\bigr)$, 
the initial form of $\varphi(\alpha^{(a)})$,
and we have 
$\,{\rm in}_\prec \bigl( \chi(\alpha^{(a)}) \bigr) \,=
\, {\rm in}_\prec \bigl( \varphi(\alpha^{(a)}) \bigr) \,$
for all $\,\alpha^{(a)} \in  {\mathcal C}^{np}_{p,m}$.
Thus image($\varphi$) and image($\chi$) have the same initial algebra, and
so we deduce the sagbi property for the polynomials $\chi(\alpha^{(a)})$
from Theorem~\ref{issagbi}.
\QED\medskip

Let ${\bf w}$ denote the weight on the variables 
${\mathcal C}^{np}_{p,m}$ defined by ${\bf w}(\alpha^{(a)}):=-a^2$.
For every incomparable pair $\gamma^{(c)},\delta^{(d)}$ in the poset 
${\mathcal C}^{np}_{p,m}$, we define the 
quadratic polynomial
$$
  R(\gamma^{(c)},\delta^{(d)})\quad:=\quad
  {\rm in}_{\bf w} \bigl( S(\gamma^{(c)},\delta^{(d)}) \bigr),
$$
where $S(\gamma^{(c)},\delta^{(d)})$ is the element of the reduced Gr\"obner
basis for the kernel of $\varphi$.
For example,  $\, R(156^{(1)},234^{(2)})\,$
equals the sum of the first ten terms in (\ref{BigSyzygy}).
The weight ${\bf w}$ is equivalent,
modulo the homogeneities of $\,{\rm kernel}(\varphi)$,
to the induced weight which was denoted by ${\mathcal A}^T \omega$ in
\cite[Chapter 11]{Sturmfels_GBCP}.  The only-if direction in 
\cite[Theorem 11.4]{Sturmfels_GBCP} implies

\begin{cor}\label{row-c-GB}
The reduced Gr\"obner basis of the kernel of $\chi$ consists 
of the quadratic polynomials $R(\gamma^{(c)},\delta^{(d)})$ as
$\gamma^{(c)},\delta^{(d)}$ run over the set of incomparable pairs in the
poset ${\mathcal C}^{np}_{p,m}$.
\end{cor}

For the Pl\"ucker ideal defining the classical Grassmannian $(n=0)$,
an explicit (but non-reduced) quadratic Gr\"obner basis is known.
It appears in the work of Hodge-Pedoe
\cite{Hodge_Pedoe} and Doubilet-Rota-Stein
 \cite{DRS}, and it  consists of the van der Waerden syzygies.
They are discussed in Gr\"obner basis language in
\cite[Section 3.1]{Sturmfels_invariant}. 
Our next aim is to introduce an analogous 
non-reduced Gr\"obner basis for the ideal
$\, {\rm kernel}(\varphi)\,$
of the quantum Grassmannian. 

We begin by defining the  skew van der Waerden syzygies
for its initial ideal 
\begin{equation}
\label{bfw}
\, {\rm kernel}(\chi)
\quad = \quad {\rm in}_{\bf w} ({\rm  kernel}(\varphi)) ,
\end{equation}
which consists of the
algebraic relations among the row-consecutive minors.
Given a sequence of integers
$\,D\,: \, 1 \leq d_1 < \cdots < d_p \leq m+p$ and any integer
$0\leq a\leq np$, let $D^{(a)}$ denote $\pm \alpha^{(a)}$, where
$\alpha$ is the reordering of the sequence $D$ and $\pm$
is the sign of the permutation which sorts the sequence $D$.
Let $T=\alpha^{(a)}\beta^{(b)}$ with $a<b$ be a  non-standard
tableau and $i$ the smallest index of a violation
$\beta_i<\alpha_{i-b+a}$.
Define increasing sequences
$$
  \begin{array}{c}
    A\quad:=\quad \alpha_1,\ldots,\alpha_{i-b+a-1}
    \qquad B\quad:=\quad \beta_{i+1},\ldots,\beta_p\\
    C\quad:=\quad \beta_1,\ldots,\beta_i,\alpha_{i-b+a},\ldots,\alpha_p.
  \end{array}
$$
For a subset $I\in\binom{[p+b-a+1]}{i}$, let $C_I$ be the corresponding
numbers from $C$ (in order) and $C_{I^c}$ be the other numbers from $C$,
also in order.
Define the {\it skew van der Waerden syzygy}
\begin{equation}\label{eq:ElQuSy}
  W(T) \quad:=\quad
  \sum_{I\in\binom{[p+b-a+1]}{i}}
  (A, C_{I^c})^{(a)}\ \cdot (C_I,B)^{(b)}\ .
\end{equation}

\begin{prop}
\label{NonRedGB}
The syzygies $\,  W(T) \,$ form a Gr\"obner basis for the kernel of $\chi$.
\end{prop}

\noindent{\bf Proof. }
Our choice of term order implies
$\,{\rm in}_\prec \bigl( W(T) \bigr) = 
T = \alpha^{(a)} \beta^{(b)}$.
Therefore it suffices to show that
$\chi(W(T))=0$.
Let $Y_1,\ldots,Y_{m+p}$ be the columns of the submatrix of 
${\mathcal L}$ given by its
rows $a+1,\ldots,b+p$. The skew van der Waerden syzygy
 $\chi(W(T))$ is an anti-symmetric, multilinear form in the
$p+b-a+1$ vectors $Y_{b_1},\ldots,Y_{b_{p+b-a+1}}$ in $(p+b-a)$-space.
\QED\medskip

The non-reduced Gr\"obner basis in Proposition \ref{NonRedGB}
can be lifted to the quantum Grassmannian
as follows. We define the {\it quantum van der Waerden syzygy}
of the non-standard tableau $T$ to 
be the unique quadratic
polynomial $V(T)$ in ${\rm kernel}(\varphi)$ which satisfies
$$ {\rm in}_{\bf w} \bigl( V(T) \bigr ) \quad = \quad
W(T), $$
and is a sum of syzygies $S(\gamma^{(c)},\delta^{(d)})$ with 
${\bf w}(\gamma^{(c)}\delta^{(d)})={\bf w}(T)$.
This syzygy exists by (\ref{bfw})
and it is unique because the quadratic generators
of the initial ideal are $k$-linearly independent,
and any two such quadratic lifts of $W(T)$ in kernel$(\varphi)$ differ by
terms whose weights are strictly less than ${\bf w}(T)$.
For instance, the quantum van der Waerden syzygy
$\, V(156^{(1)}234^{(2)})\,$ is the
polynomial with $30$ terms given in (\ref{BigSyzygy}).
It would be desirable to find an explicit formula,
perhaps in terms of the combinatorial formalism in
\cite{DRS}, for all of the skew van der Waerden syzygies $V(T)$, 
but at present we have no clue how to do this.

\vskip .1cm

The ideal of the quantum Grassmannian $K^q_{p,m}$
contains certain {\it obvious relations} which are derived
from the Grassmannian ${\it Grass}_pk^{m+p}$.
For each $\alpha\in\binom{[m+p]}{p}$ consider the polynomial
$$ g_\alpha(t) \quad = \quad
  \alpha^{(q)}\cdot t^q+ \cdots + \alpha^{(1)}\cdot t + \alpha^{(0)}\, .
$$
Given any quadratic form $F(\alpha)$ in the Pl\"ucker ideal defining
${\it Grass}_pk^{m+p}$ and any $0\leq r \leq 2q$, let $F_r$ be the
coefficient of $t^r$ in polynomial $F(g_\alpha(t))$.
Since $F(g_\alpha(t))$ is a polynomial in $t$ which vanishes identically on
$K^q_{p,m}$, each of its coefficients $F_r$ must
also vanish on $K^q_{p,m}$.
We call the collection of quadratic polynomials $F_r$ as $F$ ranges over a
generating set for the Pl\"ucker ideal of ${\it Grass}_pk^{m+p}$ the 
{\it obvious relations}. Rosenthal~\cite{Rosen94} showed the following.

\begin{prop} 
The obvious relations define $K^q_{p,m}$
set-theoretically, provided $k$ is infinite.
\end{prop}

When $p,m \leq 2$, the obvious relations coincide with
the reduced Gr\"obner basis of
Theorem~\ref{thm:gbasis}, in particular, they generate the ideal
of the quantum Grassmannian $K^q_{2,m}$.
This is no longer true for $m=p=3$.
There are 35 incomparable pairs in ${\mathcal C}^0_{3,3}$, and
hence $35$ linearly independent quadrics in the
Pl\"ucker ideal of ${\it Grass}_3k^6$.
These give rise to $35(2q+1)$ linearly independent obvious relations
but when  $q>0$ there are 
$35(2q+1)+2q-1$ incomparable pairs in ${\mathcal C}^q_{3,3}$.
Thus the obvious relations do not generate the homogeneous ideal of 
$K^q_{3,3}$.

When $q=1$ or $q=2$ then the obvious relations
generate the homogeneous ideal of $K^q_{3,3}$ together with an embedded
component supported on the irrelevant ideal.
Thus the obvious relations define $K^q_{3,3}$ scheme-theoretically, but not
ideal-theoretically. It remains an open problem whether the
the obvious relations define $K^q_{p,m}$ scheme-theoretically.

Batyrev et.al.~\cite{Batyrev} applied the familiar
sagbi property for the Grassmannian in the
construction of certain pairs of mirror 3-folds from Calabi-Yau complete
intersections in Grassmannians. We are optimistic that the
results in this paper will be similarly useful for researchers
in the fascinating interplay of algebraic geometry and theoretical physics.

The classical straightening law for the Grassmannian and its Schubert 
varieties were the starting point for the
general {\it standard monomial theory} for flag varieties.
For details and references we refer to the recent work on sagbi bases
by Gonciulea and Lakshmibai \cite{GL}.
Our results suggest 
that standard monomial theory might be
extended to certain spaces of rational curves in flag varieties
generalizing the quantum Grassmannian  $K^q_{p,m}$.

\end{document}